\documentclass{article}

\usepackage{amsmath,amsthm,amsfonts}

\newtheorem{prop}{\bf Proposition}

\newtheorem{thm}{\bf Theorem}
\newtheorem{dfn}{\bf Definition}
\newtheorem{cor}{\bf Corollary}

\newcommand{\set}[2]{\{\,#1\,|\,#2\,\}}
\newcommand{\map}[3]{#1\colon#2\rightarrow#3}
\renewcommand{\sec}[1]{\mathrm{Sec}(#1)}                %Set of sections
\newcommand{\vectorfields}[1]{{\mathfrak X}(#1)}        %Set of vectorfields
\newcommand{\Lie}[1]{{\mathcal L}_{#1}}
\newcommand{\cinfty}[1]{C^{\scriptscriptstyle\infty}(#1)}   %Set of smooth functions
\newcommand{\ver}[1]{\mathrm{Ver}(#1)}              %Set of verticals
\newcommand{\pd}[2]{\frac{\partial#1}{\partial#2}}      %Partial derivative

\newcommand{\ad}[1]{#1^\dag}                %Afine dual functor
\newcommand{\dad}[1]{\tilde{#1}}            %Dual to afine dual functor
\newcommand{\pai}[2]{\langle{#1},{#2}\rangle}       %Pairing

\def\seq 0->#1-#2->#3-#4->#5->0{0\longrightarrow
            #1\maparrow{#2}#3\maparrow{#4}#5\longrightarrow 0}
\newcommand{\maparrow}[1]{\mathrel{\mathop{\longrightarrow}\limits^{#1}}}
\def\hook{\mathop{\hbox to 6pt{\hrulefill}
                      \hbox{\vrule\phantom{\vbox to 7pt{}}}}}
\def\R{\mathit{I\kern-.35em R}}

\newcommand{\prr}{\pi}      %projection M->R
\newcommand{\pr}{\tau}      %projection E->M
\newcommand{\prv}{\bs{\pr}}     %projection V->M
\newcommand{\prad}{\ad{\pr}}    %projection Aff(E,R)->M
\newcommand{\prdad}{\dad{\pr}}  %projection Aff(E,R)^*->M
\newcommand{\prol}[2][]{\mathcal{L}^{#1}#2}     %Mixed prolongation
\newcommand{\jprol}[2][E]{\mathcal{J}^{#1}#2}   %Affine prolongation
\newcommand{\vprol}[1]{\prol[V]{#1}}            %Vector Prolongation
\newcommand{\eprol}[2][\dad{E}]{\mathcal{L}^{#1}#2} %Extended prolongation
\def\baseX{\mathcal{X}}   %Basis of prolongations
\def\baseV{\mathcal{V}}   %Basis of prolongations
\def\baseP{\mathcal{P}}   %Basis of prolongations

\newcommand{\bs}[1]{\boldsymbol{#1}}
\newcommand{\cmap}{\vartheta}

\newcommand{\sode}{pseudo-\textsc{sode}}
\newcommand{\mybox}[1]{\makebox(0,0){\footnotesize{#1}}}

\title{Lie algebroid structures and Lagrangian systems on affine bundles}

\author{Eduardo Mart\'{\i}nez$^\dag$,
 Tom Mestdag$^\ddag$ and Willy Sarlet$^\ddag$\\[10pt]
 {\small \dag Departamento de Matem\'atica Aplicada }\\
 {\small Universidad de Zaragoza, Mar\'{\i}a de Luna 3,
         50015 Zaragoza, Spain}\\[13pt]
{\small \ddag Department of Mathematical Physics and Astronomy }\\
{\small Ghent University, Krijgslaan 281, B-9000 Ghent, Belgium}}
\date{}

\begin{document}
\maketitle

\begin{abstract}
As a continuation of previous papers, we study the concept of a
Lie algebroid structure on an affine bundle by means of the
canonical immersion of the affine bundle into its bidual. We pay
particular attention to the prolongation and various lifting
procedures, and to the geometrical construction of Lagrangian-type
dynamics on an affine Lie algebroid.
\end{abstract}

\section{Introduction}

Since the book of Mackenzie \cite{Mac}, the mathematics of Lie
algebroids (and groupoids) has been studied by many authors; for a
non-exhaustive list of references, see for example
\cite{Courant,Fernandes,GrabUrbII,Tenerife1,Tenerife2,Mac2,Pop}.
The potential relevance of Lie algebroids for applications in
physics and other fields of applied mathematics has gradually
become more evident. In particular, contributions by Libermann
\cite{Liber} and Weinstein \cite{Wein} have revealed the role Lie
algebroids play in modelling certain problems in mechanics. The
concept of `Lagrangian equations' on Lie algebroids certainly
defines an interesting generalisation of Lagrangian systems as
known from classical mechanics, if only because of the more
general class of differential equations it involves while
preserving a great deal of the very rich geometrical structure of
Lagrangian (and Hamiltonian) mechanics.

One of us \cite{LMLA}, in particular, has produced evidence of this rich
structure by showing that one can prolong a Lie algebroid in such a way that
the newly obtained space has all the features of tangent bundle geometry,
which are important for the geometrical construction of Lagrangian systems.
That is to say, the prolonged Lie algebroid carries a Liouville-type section
and a vertical endomorphism which enables the definition of a
Poincar\'e-Cartan type 1-form, associated to a function $L$; the available
exterior derivative then does the rest for arriving at an analogue of the
symplectic structure from which Lagrangian equations can be derived.

In \cite{affine.algebroid}, we have started an investigation on the possible
generalisation of the concept of Lie algebroids to affine bundles. Our
principle motivation was to create a geometrical model which would be a
natural environment for a time-dependent version of Lagrange equations on Lie
algebroids, as discussed for example in \cite{Wein} and \cite{LMLA}. Since
classical time-dependent mechanics is usually described on the first-jet
bundle $J^1M$ of a manifold $M$ fibred over $\R$ (see e.g.\ \cite{Crampin,
Mangia1}), and $J^1M\rightarrow M$ is an affine bundle, it looks natural to
build up a time-dependent generalisation in such a way that $J^1M\rightarrow
M$ is the image bundle of the anchor map of a Lie algebroid structure on some
affine bundle $E\rightarrow M$. An additional indication that such a set-up is
well suited came from a naive calculus of variations approach, which gives a
clue on the analytical format one should expect for such time-dependent
Lagrange equations (see \cite{opava}).

Lie algebroids on vector bundles are known to give rise, among other things,
to a linear Poisson structure on the dual bundle, as well as a coboundary
operator on its Grassmann algebra; in fact, these properties equivalently
characterise the Lie algebroid structure. One of the features of the approach
to affine algebroids adopted in \cite{affine.algebroid} was our specific
choice to develop, in a direct way, a consistent theory of forms on sections
of an affine bundle and their exterior calculus. By contrast, however, in the
context of briefly mentioning the related Poisson structure in the concluding
remarks, we did announce a forthcoming different approach, which would be
based on the fact that an affine bundle can be regarded as an affine
sub-bundle of a vector bundle, namely the dual of its extended dual. This is
the line of reasoning we will develop here; it could be termed `indirect'
because it makes use of an imbedding into a larger bundle, but it has some
marked advantages, such as the fact that proving a number of properties
becomes much easier and that new insights come to the forefront. As a matter
of fact, one readily recognises via this approach that much (if not all) of
the theory of affine Lie algebroids can be developed without needing an extra
fibration of the base manifold over $\R$. We will accordingly start our
present analysis in this more general set-up and briefly come back to the
special case appropriate for time-dependent systems in the concluding remarks.

Note: while the very last editing of this paper was being done, we have been
informed of similar investigations on affine algebroids, which have been
carried out by Grabowski et al \cite{GGU}. The reader may find it instructive
to compare the two simultaneous developments, which are related to each other
up to and including our Section~6.

The scheme of the paper is as follows. A fairly detailed description of purely
algebraic aspects of the theory is given in Sections~2 to 4; it involves the
introduction of the concept of a Lie algebra over an affine space and aspects
of exterior calculus. A Lie algebroid structure on a general affine bundle
$\pr:E\rightarrow M$ is defined in Section~5: essentially, it comes from a
classical Lie algebroid on the dual $\dad{E}$ of the extended dual $\ad{E}$ of
$E$, with the property that the bracket of sections in the image of the
inclusion $i:E\rightarrow \dad{E}$, lies in the image of the underlying vector
bundle. Equivalent characterisations of this property can be found in the
subsequent section on the exterior differential and the associated Poisson
structure. Section~7 presents a number of simple examples of affine
algebroids. The important concept of prolongation of an algebroid is discussed
in Section~8: starting from a general construction on vector algebroids, it is
shown that the prolonged bundle inherits the affine structure coming from $E$
when the vector algebroid is the one on $\dad{E}$. For the specific case of
interest, it is further shown (see Section~9) that there is a canonical map
which gives rise to a `vertical endomorphism' on sections of the prolonged
bundle. Natural constructions which are then available are complete and
vertical lifts ; they play a role in the geometric definition of Lagrangian
systems on affine Lie algebroids, presented in Section~10.

\section{Immersion of an affine space in a vector space}

Let $A$ be an affine space modelled on a vector space $V$, and let
$\ad{A}=\mathrm{Aff}(A,\R)$ be the extended dual of $A$, that is, the vector
space of all affine maps from $A$ to the real line. We consider the bidual
$\dad{A}$ of $A$, in the sense $\dad{A}=(\ad{A})^*$. It is well known that in
the case of a vector space $V$ , the bidual $\dad{V}=(V^*)^*$ is isomorphic to
$V$. In the case of an affine space, the bidual includes `a copy' of $A$, as
it is shown in the following statement.

\begin{prop}
The map $\map{i}{A}{\dad{A}}$ given by $i(a)(\varphi)=\varphi(a)$ is an
injective affine map, whose associated vector map is
$\map{\bs{i}}{V}{\dad{A}}$ given by
$\bs{i}(\bs{v})(\varphi)=\bs{\varphi}(\bs{v})$
\end{prop}
\begin{proof}
If $a\in A$ and $\bs{v}\in V$, then for all $\varphi\in\ad{A}$,
$$
i(a+\bs{v})(\varphi)=\varphi(a+\bs{v})=\varphi(a)+\bs{\varphi}(\bs{v})
=i(a)(\varphi)+\bs{i}(\bs{v}),
$$
from which it follows that $i$ is an affine map whose associated linear map is
$\bs{i}$.

To prove that $i$ is injective, it suffices to prove that $\bs{i}$ is
injective, which is obvious since if $\bs{v}$ is an element in the kernel of
$\bs{i}$ then $\bs{i}(\bs{v})(\bs{\varphi})=\bs{\varphi}(\bs{v})=0$ for all
$\bs{\varphi}\in V^*$, hence $\bs{v}=0$.
\end{proof}

The vector space $\dad{A}$ is foliated by hyperplanes parallel to the image of
$\bs{i}$. Every vector $z\in\dad A$ is either of the form $z=\bs{i}(\bs{v})$
for some $\bs{v}\in V$ or of the form $z=\lambda i(a)$ for some
$\lambda\in\R\setminus\{0\}$ and $a\in A$. Moreover, $\lambda, a$ and $\bs{v}$
are uniquely determined by $z$. The image of the map $i$ consists of the
points for which $\lambda=1$. To understand this description in more detail,
we will prove that we have an exact sequence of vector spaces $\seq 0-> V
-\bs{i}-> \dad{A} - -> \R->0$. To this end we consider the dual sequence.

\begin{prop}
Let $\map{l}{\R}{\ad{A}}$ be the map that associates to $\lambda\in\R$ the
constant function $\lambda$ on $A$. Let $\map{k}{\ad{A}}{V^*}$ the map that
associates to every affine function on $A$ the corresponding linear function on
$V$. Then, the sequence of vector spaces
$$
\seq 0-> \R -l-> \ad{A} -k-> V^* ->0
$$
is exact.
\end{prop}
\begin{proof}
Indeed, it is clear that $l$ is injective, $k$ is surjective and $k\circ l=0$,
so that $\mathrm{Im}(l)\subset\mathrm{Ker}(k)$. If $\varphi\in\ad{A}$ is in the
kernel of $k$, that is the linear part of $\varphi$ vanishes, then for every
pair of points $a$ and $b=a+\bs{v}$ we have that
$$
\varphi(b)=\varphi(a+\bs{v})=\varphi(a)+\bs{\varphi}(\bs{v})=\varphi(a),
$$
that is $\varphi$ is constant, and hence in the image of $l$.
\end{proof}

The dual map of $k$ is $\bs{i}$, since for $\bs{v}\in V$ we have
$$
\pai{k(\varphi)}{\bs{v}}
  =\pai{\bs{\varphi}}{\bs{v}}
  =\pai{\varphi}{\bs{i}(\bs{v})}
$$
The dual map $j$ of the map $l$ is given by $j(\alpha
i(a)+\bs{i}\bs{v})=\alpha$. Indeed, for every $\lambda\in\R$ we have
$$
j(z)\lambda
  =\pai{z}{l(\lambda)}
  =\pai{\alpha i(a)+\bs{i}(\bs{v})}{l(\lambda)}
  =\alpha\pai{i(a))}{l(\lambda)} +\pai{\bs{i}(\bs{v})}{l(\lambda)}
  =\alpha\lambda
$$
It follows that

\begin{cor}
If $A$ is finite dimensional, then the sequence
$$
\seq 0-> V -\bs{i}-> \dad{A} -j-> \R->0
$$
is exact.
\end{cor}

Note that in this way we can clearly identify the image of $V$ as the
hyperplane of $\dad{A}$ with equation $j(z)=0$, and the image of $A$ as the
hyperplane of $\dad{A}$ with equation $j(z)=1$, in other words
$$
\bs{i}(V)=j^{-1}(0)\qquad\text{and}\qquad i(V)=j^{-1}(1).
$$
Note in passing that if we have an exact sequence $\seq 0-> V -\alpha-> W -j->
\R->0$, then we can define $A=j^{-1}(1)$; it follows that $A$ is an affine
space modelled on the vector space $V$ and $W$ is canonically isomorphic to
$\dad{A}$. The isomorphism is the dual map of $\Psi:W^*\rightarrow\ad{A}, \
\Psi(\phi)(a)=\phi(i(a))$, where $i:A\rightarrow W$ is the canonical
inclusion.

\medskip

We now discuss the construction of a basis for $\ad{A}$. Let
$(O,\{\bs{e}_i\})$ be an affine frame on $A$. Thus every point $a$ has a
representation $a=O+v^i\bs{e}_i$. The family of affine maps
$\{e^0,e^1,\ldots,e^n\}$ given by
$$
e^0(a)=1 \qquad\qquad e^i(a)=v^i,
$$
is a basis for $\ad{A}$. If $\varphi\in A^ \dag$, and we put
$\varphi_0=\varphi(O)$ and $\varphi_i=\bs{\varphi}(\bs{e}_i)$, then
$\varphi=\varphi_0e^0+\varphi_ie^i$.

It is to be noticed that, contrary to $e^1,\ldots,e^n$, the map $e^0$ does not
depend on the frame we have chosen for $A$. In fact, $e^0$ coincides with the
map $j$.

Let now $\{e_0,e_1,\ldots,e_n\}$ denote the basis of $\dad{A}$ dual
to $\{e^0,e^1,\ldots,e^n\}$. Then the image of the canonical
immersion is given by
$$
i(O)=e_0\qquad\qquad \bs{i}(\bs{e}_i)=e_i
$$
from which it follows that for $a=O+v^i\bs{e}_i$, we have $i(a)=e_0+v^ie_i$.
If we denote by $(x^0,x^1,\ldots,x^n)$ the coordinate system on
$\dad{A}$ associated to the basis $\{e_0,\ldots,e_n\}$, then the equation of the
image of the map $i$ is $x^0=1$, while the equation of the image of $\bs{i}$ is
$x^0=0$.

Coordinates in $\dad{A}^*=\ad{A}$ associated to the basis above will be
denoted by $(\mu_0,\mu_1,\ldots,\mu_n)$, that is
$\mu_\alpha(\varphi)=\pai{\bs{e}_\alpha}{\varphi}$ for every
$\varphi\in\ad{A}$.

\section{Lie algebra structure over an affine space}

Implicit in our previous paper~\cite{affine.algebroid} is the following
definition of a Lie algebra over an affine space.

\begin{dfn}
Let $A$ be an affine space over a vector space $V$. A Lie algebra structure on
$A$ is given by
\begin{itemize}
\item a Lie algebra structure $[\,,\,]$ on $V$, and
\item an action by derivations of $A$ on $V$, i.e.\ a map $\map{D}{A\times
V}{V}$, $(a,\bs{v})\mapsto D_a\bs{v}$ with the properties
\begin{align*}
&D_a(\lambda \bs{v})=\lambda D_a\bs{v}\\
&D_a(\bs{v}+\bs{w})=D_a\bs{v}+D_a\bs{w}\\
&D_a[\bs{v},\bs{w}]=[D_a\bs{v},\bs{w}]+[\bs{v},D_a\bs{w}],
\end{align*}
\item satisfying the compatibility property
$$
D_{a+\bs{v}}\bs{w}=D_a\bs{w}+[\bs{v},\bs{w}].
$$
\end{itemize}
\end{dfn}
Incidentally, it is sufficient to require in the first item that the bracket
on $V$ is $\R$-bilinear and skew-symmetric, since the Jacobi identity then
follows from the requirements on $D_a$.

If we use a bracket notation $[a,\bs{v}]\equiv D_a\bs{v}$, then the conditions
in the definition above read
\begin{align*}
&[a,\lambda \bs{v}]=\lambda [a,\bs{v}]   \\
&[a,\bs{v}+\bs{w}]=[a,\bs{v}]+[a,\bs{w}]\\
&[a,[\bs{v},\bs{w}]]=[[a,\bs{v}],\bs{w}]+[\bs{v},[a,\bs{w}]]\\
&[a+\bs{v},\bs{w}]=[a,\bs{w}]+[\bs{v},\bs{w}].
\end{align*}
This allows us to define a bracket of elements of $A$ by putting
$$
[\bs{v},a]=-[a,\bs{v}]
$$
and, if $b=a+\bs{v}$, then
$$
[a,b]=[a,\bs{v}].
$$
This bracket is skew-symmetric by construction and also satisfies a
Jacobi-type property.

\begin{thm}
A Lie algebra structure over an affine space $A$ is equivalent to a Lie
algebra extension of the trivial Lie algebra $\R$ by $V$. Explicitly, it is
equivalent to the exact sequence of vector spaces $\seq 0-> V -i-> \dad{A}
-j-> \R->0$ being an exact sequence of Lie algebras.
\end{thm}
\begin{proof}
If the exact sequence is one of Lie algebras, we of course have a Lie algebra
structure on $V$ and the map $D$ determined by
$D_a\bs{v}=[i(a),\bs{i}(\bs{v})]$ satisfies all requirements to define a Lie
algebra structure on $A$.

Conversely, assume we have a Lie algebra structure on the affine space $A$. If
we fix an element $a\in A$, then every element $z\in\dad{A}$ can be written in
the form $z=\lambda i(a)+\bs{i}(\bs{v})$. We can define a bracket of two
elements $z_1=\lambda_1 i(a)+\bs{i}(\bs{v}_1)$ and $z_2=\lambda_2
i(a)+\bs{i}(\bs{v}_2)$ by
$$
[z_1,z_2]
 =\bs{i}([\bs{v}_1,\bs{v}_2]+\lambda_1 D_a\bs{v}_2-\lambda_2 D_a\bs{v}_1)
$$
This bracket is clearly bi-linear and skew-symmetric, and a straightforward
calculation shows that it satisfies the Jacobi identity. Moreover, the
definition does not depend on the choice of the point $a$; if $a'$ is another
point in $A$, then $a'=a+\bs{w}$ for some $\bs{w}\in V$ and the compatibility
condition implies that the result is independent of that choice. Finally, it
is obvious that the maps $i$ and $j$ then are Lie algebra homomorphisms.
\end{proof}

Notice that the only condition for a Lie algebra structure on $\dad{A}$ to be
an extension of $\R$ by $V$ is that the bracket takes values in $V$,
symbolically: $[\dad{A},\dad{A}]\subset V$.

Once we have chosen an affine frame on $A$, we have that the bracket on
$\dad{A}$ is determined by the brackets of the associated basis elements.
These must be of the form
$$
[e_0,e_0]=0\qquad [e_0,e_j]= C^k_{0j}e_k\qquad [e_i,e_j]=C^k_{ij}e_k,
$$
since all brackets must take values in the image of the map $\bs{i}$.

\medskip

It is well known that a Lie algebra structure on a vector space defines, and
is defined by, a linear Poisson structure on the dual vector space. In the
light of the results of the previous section we have a Poisson bracket on
$V^*$ and one on $\ad{A}$. Furthermore, $\dad{A}$ being an extension of $\R$
by $V$, we have that the Poisson structure $\Lambda_{\ad{A}}$ is an extension
by $\Lambda_{V^*}$ of $\Lambda_{\R}=0$ (see~\cite{poisson.extension} for the
details on Poisson extensions and their relations to Lie algebra extensions).
Therefore, once we have fixed a point $a\in A$, we have a splitting of the
sequence $\seq 0->V-i->\dad{A}-j->\R->0$ given by $h(\lambda)=\lambda i(a)$,
and the Poisson tensor can be written as $\Lambda_{\ad{A}}=\Lambda_{V^*}+
X_a\wedge X_{D_a}$, where $X_{D_a}$ is the linear vector field associated to
the linear map $D_a\in\mathrm{End}(V)$, and $X_a$ is the constant vector
corresponding to $i(a)$. Moreover, we have that
$\mathcal{L}_{X_{D_a}}\Lambda_{V^*}=0$.

In the coordinates $\mu_0,\mu_1,\ldots,\mu_n$ on $\ad{A}$ associated to the
basis $\{e^\alpha\}$, we have
$$
\{\mu_0,\mu_0\}=0 \qquad \{\mu_0,\mu_j\}=\mu_k\,C^k_{0j} \qquad
\{\mu_i,\mu_j\}=\mu_k\,C^k_{ij},
$$
where $C^k_{0j}, C^k_{ij}$ are the structure constants introduced above.
Therefore, the Poisson tensor reads
$$
\Lambda_\ad{A}=\frac{1}{2}\mu_k\,C^k_{ij}\pd{}{\mu_i}\wedge \pd{}{\mu_j}
              +\mu_k\,C^k_{0j}\pd{}{\mu_0}\wedge \pd{}{\mu_j}\,,
$$
where the first term is $\Lambda_V$ and the second one is  $X_a\wedge X_{D_a}$
with $a=O$, the origin.

\section{Exterior algebra over an affine space}

In~\cite{affine.algebroid} we have defined the concept of differential forms
on sections of an affine bundle. We aim in this section to establish at the
algebraic level, the relation with ordinary forms on the vector extension. We
will re-state in the algebraic setting the definition of a $k$-form, $k\geq1$.

\begin{dfn}
A $k$-form on an affine space $A$ is a map $\map{\omega}{A\times\cdots\times
A}{\R}$ for which there exists a $k$-form $\bs{\omega}$ on the associated
vector space $V$, and a map $\map{\omega_0}{A\times V\times\cdots\times
V}{\R}$ with the following properties
\begin{enumerate}
\item $\omega_0$  is skew-symmetric and linear in its $k-1$ vector arguments.
\item For every $a\in A$ and for every $\bs{v}_1,\bs{v}_2,\ldots,\bs{v}_{k}\in
V$, we have
$$
\omega_0(a+\bs{v}_1,\bs{v}_2,\ldots,\bs{v}_{k})
=\omega_0(a,\bs{v}_2,\ldots,\bs{v}_{k})+
\bs{\omega}(\bs{v}_1,\bs{v}_2,\ldots,\bs{v}_{k}).
$$
\item For every $a_1,\ldots,a_k\in A$, if we choose an arbitrary $a_0\in A$ and
put $a_i=a_0+\bs{v}_i$, we have
$$
\omega(a_1,\ldots,a_{k}) =\sum_{j=1}^k(-1)^{j+1}
\omega_0(a_0,\bs{v}_1,\ldots,\widehat{\bs{v}}_j,\ldots,\bs{v}_{k})
+ \bs{\omega}(\bs{v}_1,\bs{v}_2,\ldots,\bs{v}_{k}).
$$
\end{enumerate}
\end{dfn}

We next show that a $k$-form on $A$ is just the pull-back by the canonical
immersion of a $k$-form on $\dad{A}$, in other words, $\omega$ is a $k$-form
if we have
$$
\omega(a_1,\ldots,a_k)=\dad{\omega}(i(a_1),\ldots,i(a_k)),
$$
for some ordinary exterior $k$-form $\dad{\omega}$ on the vector space
$\dad{A}$.

\begin{prop}
If $\dad{\omega}$ is a $k$-form on $\dad{A}$ then $\omega=i^*\dad{\omega}$ is
a $k$-form on the affine space $A$. Conversely,  given a $k$-form on the
affine space $A$, there exists a unique $k$-form $\dad{\omega}$ on $\dad{A}$
such that $\omega=i^*\dad{\omega}$.
\end{prop}
\begin{proof}
For a given $k$-form $\dad{\omega}$ on $\dad{A}$ we define the maps
\begin{align*}
&\omega(a_1,\ldots,a_k)=\dad{\omega}(i(a_1),\ldots,i(a_k)),\\
&\omega_0(a,\bs{v}_2,\ldots,\bs{v}_k)
  =\dad{\omega}(i(a),\bs{i}(\bs{v}_2)\ldots,\bs{i}(\bs{v}_k)),\\
&\bs{\omega}(\bs{v}_1\ldots,\bs{v}_k)
  =\dad{\omega}(\bs{i}(\bs{v}_1)\ldots,\bs{i}(\bs{v}_k)).
\end{align*}
Then, conditions \textit{1} and \textit{2} in the definition above are
trivially satisfied. Moreover, if we fix $a_0\in A$ and write
$a_i=a_0+\bs{v}_i$, then by skew-symmetry of $\dad{\omega}$ we have
\begin{align*}
\omega(a_1,\ldots,a_k)
&=\dad{\omega}(i(a_1),\ldots,i(a_k))\\
&=\dad{\omega}(i(a_0)+\bs{i}(\bs{v}_1),\ldots,i(a_0)+\bs{i}(\bs{v}_k))\\
&=\sum_{j=1}^k(-1)^{j+1}
  \dad{\omega}(i(a_0),\bs{i}(\bs{v}_1),\ldots,
  \widehat{\bs{i}(\bs{v}_j)},\ldots,\bs{i}(\bs{v}_k))
+\dad{\omega}(\bs{i}(\bs{v}_1),\ldots,\bs{i}(\bs{v}_k))\\
&=\sum_{j=1}^k(-1)^{j+1}
\omega_0(a_0,\bs{v}_1,\ldots,\widehat{\bs{v}_j},\ldots,\bs{v}_k)
+\bs{\omega}(\bs{v}_1,\ldots,\bs{v}_k),
\end{align*}
which proves condition \textit{3}.

Conversely, assume we are given a $k$-form $\omega$ on the affine space $A$
with its associated $\omega_0$ and $\bs{\omega}$. Fixing $a_0\in A$, we know
that every point $z\in\dad{A}$ can be written in the form $z=\lambda
i(a_0)+\bs{i}(\bs{v})$ for $\lambda\in\R$ and $\bs{v}\in V$. We define the map
$\dad{\omega}$ by
\begin{align*}
\dad{\omega}(z_1,\ldots,z_k)
&=\dad{\omega}(\lambda_1 i(a_0)+\bs{i}(\bs{v}_1),\ldots,
                \lambda_k i(a_0)+\bs{i}(\bs{v}_k))\\
&=\sum_{j=1}^k(-1)^{j+1}\lambda_j\,
\omega_0(a_0,\bs{v}_1,\ldots,\widehat{\bs{v}_j},\ldots,\bs{v}_k)
+\bs{\omega}(\bs{v}_1,\ldots,\bs{v}_k).
\end{align*}
By virtue of \textit{1} and \textit{3}, it follows that $\dad{\omega}$ is
multi-linear and skew-symmetric, i.e.\ it is a $k$-form on $\dad{A}$.
Moreover, $\dad{\omega}(z_1,\ldots,z_k)$ is independent of the choice of the
point $a_0$. Indeed, if we choose a different point $a'_0=a_0+\bs{w}$, then
$z_j=\lambda_j i(a'_0)+\bs{v}'_j$ with $\bs{v}'_j=\bs{v}_j-\lambda_j\bs{w}$,
and applying the definition above we get
\begin{align*}
\dad{\omega}(z_1,\ldots,z_k)
&=\dad{\omega}(\lambda_1 i(a'_0)+\bs{i}(\bs{v}'_1),\ldots,
                 \lambda_k i(a'_0)+\bs{i}(\bs{v}'_k))\\
&=\sum_{j=1}^k(-1)^{j+1}\lambda_j\,
\omega_0(a'_0,\bs{v}'_1,\ldots,\widehat{\bs{v}'_j},\ldots,\bs{v}'_k)
+\bs{\omega}(\bs{v}'_1,\ldots,\bs{v}'_k)\\
&=\sum_{j=1}^k(-1)^{j+1}\lambda_j\, \omega_0(a_0+\bs{w},\bs{v}_1-
\lambda_1\bs{w},\ldots,\widehat{\bs{v}_j},\ldots,\bs{v}_k-\lambda_k\bs{w})\\
&\hspace*{1cm} +\bs{\omega}(\bs{v}_1-\lambda_1\bs{w},\ldots,\bs{v}_k-\lambda_k\bs{w}) \\
&=\sum_{j=1}^k(-1)^{j+1}\lambda_j\,
\omega_0(a_0,\bs{v}_1,\ldots,\widehat{\bs{v}_j},\ldots,\bs{v}_k)
+\bs{\omega}(\bs{v}_1,\ldots,\bs{v}_k)
\end{align*}
where we have used the properties of $\omega_0$ and $\bs{\omega}$.

The form $\dad{\omega}$ is unique, since, if $\dad{\theta}$ is a $k$-form on
$\dad{A}$ such that $i^*\dad{\theta}=0$, then it follows that the associated
$\theta_0$ and $\bs{\theta}$ vanish from where we deduce that
$\dad{\theta}=0$.
\end{proof}

Once a reference frame has been fixed on $A$, a 1-form $\dad{\omega}$ on
$\dad{A}$ is of the form $\dad{\omega}=\omega_0e^0+\omega_ie^i$, and then the
local representation of $\omega=i^*\dad\omega$ looks exactly the same.

More generally,  a p-form $\dad{\omega}$ on $\dad{A}$ is of the form
$$
\dad{\omega}=\frac{1}{p!}\sum_{\mu_1,\ldots,\mu_p=0}^n
\omega_{\mu_1\cdots\mu_p}e^{\mu_1}\wedge\cdots\wedge e^{\mu_p}
$$
and thus
$$
\omega=\frac{1}{(p-1)!}\,\omega_{0i_1\cdots i_{p-1}}e^0\wedge
e^{i_1}\wedge\cdots\wedge e^{i_{p-1}} +\frac{1}{p!}\,\omega_{i_1\cdots
i_p}e^{i_1}\wedge\cdots\wedge e^{i_p},
$$
where the first term corresponds to the sum of terms involving $\omega_0$ and
the second to $\bs{\omega}$.

\section{Lie algebroid structures over affine bundles}

Let $\map{\pr}{E}{M}$ be an affine bundle with associated vector bundle
$\map{\prv}{V}{M}$. We consider the bundle $\map{\prad}{\ad{E}}{M}$ whose
fibre over $m\in M$ is the extended dual $\ad{E}_m$ of the fibre $E_m$. We
also consider the dual bundle $\map{\prdad}{\dad{E}}{M}$, whose fibre at $m$
is $\dad{E}_m=(\ad{E}_m)^*$. At every point $m$ we have the exact sequence of
vector spaces
$$
\seq 0-> V_m -\bs{i}_m-> \dad{E}_m -j_m-> \R ->0
$$
and therefore an exact sequence of vector bundles over $M$
$$
\seq 0-> V -\bs{i}-> \dad{E} -j-> M\times\R ->0.
$$
On the other hand there is the canonical immersion $\map{i}{E}{\dad{E}}$, so
that $i(E_m)=j^{-1}((m,1))$.

By taking sections, we have the exact sequence of real vector spaces (and
$\cinfty{M}$-modules)
$$
\seq 0-> \sec{V} -\bs{i}-> \sec{\dad{E}} -j-> \cinfty{M} ->0.
$$
and an inclusion $\map{i}{\sec{E}}{\sec{\dad{E}}}$, whereby we make no
notational distinction between the bundle maps and the induced maps of
sections (i.e.\ if $\sigma$ is a section and $r$ is a bundle map over the
identity, we write $r(\sigma)$ instead of $r\circ\sigma$.). It follows that if
we fix a section $\sigma$ of $E$ then we have a splitting of the above
sequence and any section $\zeta$ of $\dad{E}$ can be written as
$\zeta=fi(\sigma)+\bs{i}(\bs{\eta})$, for some section $\bs{\eta}$ of $V$ and
where $f=j(\zeta)$.

Also, since $V$ is the vector bundle associated to the affine bundle $E$, we
have that $\sec{V}$ is the (real) vector space associated to the affine space
$\sec{E}$.

\medskip

In~\cite{affine.algebroid} we defined the concept of Lie algebroid structure
on the class of affine bundles whose base manifold is further fibred over
$\R$. Here we start with a more general definition, which in the light of the
previous section can be expressed as follows.

\begin{dfn}
A Lie algebroid structure on $E$ consists of a Lie algebra structure on the
(real) affine space of sections of $E$ together with an affine map
$\map{\rho}{E}{TM}$ (the anchor), satisfying the following compatibility
condition
$$
D_\sigma(f\bs{\zeta})
       =\rho(\sigma)(f)\,\bs{\zeta}+fD_\sigma\bs{\zeta},
$$
for every $\sigma\in\sec{E}$, $\bs{\zeta}\in\sec{V}$ and $f\in\cinfty{M}$, and
where $D_\sigma$ is the action $\sigma\mapsto D_\sigma$ of $\sec{E}$ on
$\sec{V}$.
\end{dfn}

The compatibility condition ensures that the association $\sigma\mapsto
D_\sigma$, which acts by derivations on the real Lie algebra $\sec{V}$, also
acts by derivations on the $\cinfty{M}$-module $\sec{V}$.

The anchor map $\rho$ extends to a linear map $\dad{\rho}: \dad{E}\rightarrow
TM$, which we will describe in more detail below. It is of interest, however,
to observe now already that the map $\map{\bs{i}}{V}{\dad{E}}$ is a morphism
of Lie algebroids, since we have
$$
[\bs{i}(\bs{\eta}_1),\bs{i}(\bs{\eta}_2)]=\bs{i}([\bs{\eta}_1,\bs{\eta}_2])
\quad\text{and}\quad \dad{\rho}\circ \bs{i}=\bs{\rho},
$$
where $\bs{\rho}$ is the linear part of $\rho$. On the contrary, the map
$\map{j}{\dad{E}}{M\times\R}$ is \textsc{not} a morphism of Lie algebroids,
since we have that
$j([\zeta_1,\zeta_2])=\dad{\rho}(\zeta_1)f_2-\dad{\rho}(\zeta_2)f_1$, while
$[j(\zeta_1),j(\zeta_2)]=0$ since the fibres of $M\times\R$ are 1-dimensional.

The affine Lie algebroid structure we studied in~\cite{affine.algebroid} is
the case that $M$ is further fibred over the real line $\map{\prr}{M}{\R}$ and
the anchor map $\rho$ takes values in $J^1M$. Notice that such an extra
fibration is not generally available, not even locally. For instance, if we
take any affine bundle $\map{\pr}{E}{M}$ with the trivial Lie algebroid
structure (null bracket and anchor) then there is no fibration over $\R$ such
that the image of $\rho$ is in $J^1M$, since the vectors in $i(J^1M)$ are
non-zero.

The following result shows that one can alternatively define an affine Lie
algebroid structure on $E$ as a vector Lie algebroid structure
$([\,,\,],\dad{\rho})$ on $\dad{E}$ such that the bracket of two sections in
the image of $i$ belongs to the image of $\bs{i}$.

\begin{thm}
A Lie algebroid structure on the vector bundle $\map{\prdad}{\dad{E}}{M}$
which is such that the bracket of sections in the image of $i$ lies in the
image of $\bs{i}$ defines a Lie algebroid structure on the affine bundle
$\map{\pr}{E}{M}$, whereby the brackets and maps are determined by the
following relations:
\begin{gather*}
\bs{i}([\bs{\eta}_1,\bs{\eta}_2])=[\bs{i}(\bs{\eta}_1),\bs{i}(\bs{\eta}_2)]\\
\bs{i}(D_\sigma\bs{\eta})=[i(\sigma),\bs{i}(\bs{\eta})]\\
\rho(\sigma)=\dad{\rho}(i(\sigma)).
\end{gather*}

Conversely, a Lie algebroid structure on the affine bundle $\map{\pr}{E}{M}$
extends to a Lie algebroid structure on the vector bundle
$\map{\prdad}{\dad{E}}{M}$ such that the bracket of sections in the image of
$i$ is in the image of $\bs{i}$. If we fix a section $\sigma$ of $E$ and write
sections $\zeta$ of $\dad{E}$ (locally) in the form
$\zeta=fi(\sigma)+\bs{i}(\bs{\eta})$ then the anchor and the bracket are given
by
\begin{gather*}
\dad{\rho}(\zeta)=f\rho(\sigma)+\bs{\rho}(\bs{\eta})\\
[\zeta_1,\zeta_2]=\Bigl(\dad{\rho}(\zeta_1)(f_2)-\dad{\rho}(\zeta_2)(f_1)\Bigr)i(\sigma)
+\bs{i}\Bigl([\bs{\eta}_1,\bs{\eta_2}]+f_1D_\sigma\bs{\eta}_2-f_2D_\sigma\bs{\eta}_1\Bigr).
\end{gather*}
\end{thm}
\begin{proof} The verification of the above statements is straightforward but
rather lengthy. We limit ourselves to checking that the compatibility
conditions between brackets and anchors are satisfied. For the first part, we
find
\begin{align*}
\bs{i}\Bigl(D_\sigma(f\bs{\eta})\Bigr)&=[i(\sigma),\bs{i}(f\bs{\eta})]
=[i(\sigma),f\bs{i}(\bs{\eta})]\\
&=\dad{\rho}(i(\sigma))(f)\bs{i}(\bs{\eta})+f[i(\sigma),\bs{i}(\bs{\eta})]
=\bs{i}\Bigl(\rho(\sigma)(f)\bs{\eta}+fD_\sigma(\bs{\eta})\Bigr),
\end{align*}
form which it follows that $D_\sigma(f\bs{\eta})=\rho(\sigma)(f)\bs{\eta} +
fD_\sigma(\bs{\eta})$.

For the converse, observe that
\begin{align*}
[\zeta_1,f\zeta_2]-f[\zeta_1,\zeta_2]
&=f_2\dad{\rho}(\zeta_1)(f)i(\sigma)+\bs{i}\Big(f_1\rho(\sigma)(f)\bs{\eta_2}
  +\bs{\rho}(\bs{\eta_1})(f)\bs{\eta}_2\Big)\\
&=\dad{\rho}(\zeta_1)(f)f_2i(\sigma)+(f_1\rho(\sigma) + \bs{\rho}(\bs{\eta_1}))(f)\bs{i}(\bs{\eta}_2)\\
&=\dad{\rho}(\zeta_1)(f)\zeta_2,
\end{align*}
which is the required compatibility condition.
\end{proof}

In coordinates, if $x^i$ denote coordinates on $M$ and $y^\alpha$ fibre
coordinates on $E$ with respect to some local frame $(e_0,\{\bs{e}_\alpha\})$
of sections of $E$, then we have
\begin{align*}
\rho(e_0+y^\alpha \bs{e}_\alpha)&=(\rho^i_0+\rho^i_\alpha
y^\alpha)\pd{}{x^i}\,,\\
[\bs{e}_\alpha,\bs{e}_\beta]&=C^\gamma_{\alpha\beta}\bs{e}_\gamma\,,\\
D_{e_0}\bs{e}_\beta&=C^\gamma_{0\beta}\bs{e}_\gamma,
\end{align*}
for some functions $\rho^i_0$, $\rho^i_\alpha$, $C^\gamma_{0\beta}$ and
$C^\gamma_{\alpha\beta}$ on $M$.

Taking the local basis of sections of $\dad{E}$, associated to the above
frame, it follows that
$$
\dad{\rho}(y^0e_0+y^\alpha e_\alpha)=(\rho^i_0y^0+\rho^i_\alpha
y^\alpha)\pd{}{x^i}\,,
$$
and the bracket is determined by
$$
[e_0,e_0]=0
\qquad
[e_0,e_\beta]= C^\gamma_{0\beta}e_\gamma
\qquad
[e_\alpha,e_\beta]=C^\gamma_{\alpha\beta}e_\gamma.
$$

As a final remark we mention that the orbits of the Lie algebroid
$\map{\bs{\pr}}{V}{M}$ are subsets of  the orbits of the algebroid
$\map{\dad{\pr}}{\dad{E}}{M}$ and they are equal if and only if there exists a
section $\sigma$ of $E$ such that $\rho(\sigma)$ is in the image of
$\bs{\rho}$.

\section{Exterior differential and Poisson structure}

Now that we have proved that a Lie algebroid structure on an affine bundle is
equivalent to a Lie algebroid structure on $\dad{E}$, we can define the
exterior differential operator on $E$ by pulling back the exterior
differential on $\dad{E}$. More precisely, given a $k$-form $\omega$ on the
affine bundle $E$ we know that there exists a unique $\dad{\omega}$ on
$\dad{E}$ such that $\omega=i^*(\dad{\omega})$. Then we define $d\omega$ as
the $(k+1)$-form given by
$$
d\omega=i^*(d\dad{\omega}).
$$
It is easy to see that this definition is equivalent to the one given
in~\cite{affine.algebroid} (at least when $M$ is fibred over $\R$).

The definition given here has some clear advantages. For instance, the
property $d^2=0$ which was rather difficult to prove
in~\cite{affine.algebroid}, becomes evident now:
$$
d^2\omega =d(d(i^*\dad{\omega})) =d(i^*d\dad{\omega})
 =i^*(d^2\dad{\omega})=0.
$$

The differential $\bs{d}$ on the Lie algebroid $V$ is also related to the
differential on $\dad{E}$; for every $k$-form $\dad{\omega}$ on $E$ we have
that
$$
\bs{d}\bs{i}^*\dad{\omega}=\bs{i}^*(d\dad{\omega}),
$$
which in fact simply expresses that $\bs{i}$ is a morphism of Lie algebroids.
Hence, one can find the differential of a form $\bs{\omega}$ on $V$ by
choosing a $k$-form $\dad{\omega}$ on $\dad{E}$ such that
$\bs{i}^*\dad{\omega}=\bs{\omega}$ and then obtain $\bs{d}\bs{\omega}$ as
$\bs{i}^*(d\dad{\omega})$. That this does not depend on the choice of
$\dad{\omega}$ is expressed in equivalent terms by the following result.

\begin{prop}
The ideal $ \mathcal{I}=\{\dad{\omega}\,|\,\bs{i}^*\dad{\omega}=0\}$ is a
differential ideal, i.e.\ $d\mathcal{I}\subset\mathcal{I}$.

\end{prop}
\begin{proof}
If $\bs{i}^*\dad{\omega}=0$, then the same is true for $d\dad{\omega}$ since
$\bs{i}^*(d\dad{\omega})=\bs{d}\bs{i}^*(\dad{\omega})=0$.
\end{proof}

The ideal $\mathcal{I}$ is generated by the 1-form $e^0$,
$\mathcal{I}=\{e^0\wedge\theta\mid\text{$\theta$ is a form on $\dad{E}$}\}$ so
that $de^0$ belongs to $\mathcal{I}$. The following result shows that $e^0$ in
fact is $d$-closed and that this property characterises affine structures.

\begin{thm}
A Lie algebroid structure on $\dad{E}$ restricts to a Lie algebroid structure
on the affine bundle $E$ if and only if the exterior differential satisfies
$de^0=0$.
\end{thm}
\begin{proof}
Indeed, taking two sections $\sigma_1$ and $\sigma_2$ of $E$ we have
\begin{align*}
de^0(i(\sigma_1),i(\sigma_2)) &=\dad{\rho}(\sigma_1)\pai{e^0}{i(\sigma_2)}
-\dad{\rho}(\sigma_2)\pai{e^0}{i(\sigma_1)}
-\pai{e^0}{[i(\sigma_1),i(\sigma_2)]}\\
&=-\pai{e^0}{[i(\sigma_1),i(\sigma_2)]}
\end{align*}
It follows that $[i(\sigma_1),i(\sigma_2)]$ is in
$\textrm{Im}(\bs{i})=\mathrm{Ker}(e^0)$ if and only if $de^0$ vanishes on the
image of $i$, which spans $\dad{E}$.
\end{proof}

In coordinates, the exterior differential operator is determined by
$$
df=\rho^i_0\pd{f}{x^i}e^0+\rho^i_\alpha\pd{f}{x^i}e^\alpha, \qquad \mbox{for\
\ } f\in\cinfty{M}
$$
and
$$
de^0=0,\qquad de^\gamma=-C^\gamma_{0\beta}\,e^0\wedge e^\beta
-\frac{1}{2}C^\gamma_{\alpha\beta}\,e^\alpha\wedge e^\beta.
$$

\bigskip

In the special case that $\rho(E)\subset J^1M$ ($M$ fibred over $\R$), we have
that $dt=e^0$ so that $e^0$ is not only closed but also exact. In fact, this
is the condition for a Lie algebroid structure on an affine bundle to have a
1-jet-valued anchor. Indeed, if there exists an $f\in\cinfty{M}$ such that
$df=e^0$, then the partial derivatives of $f$ cannot simultaneously vanish,
hence $f$ defines a local fibration and then for any section $\sigma$ of $E$
we have that $\rho(\sigma)f=\pai{df}{\sigma}=\pai{e^0}{\sigma}=1$, which is
the condition for the anchor being 1-jet-valued.

\bigskip

When we have a Lie algebroid structure on $\dad{E}$, there is a Poisson
bracket on the dual bundle $\dad{E}^*=\ad{E}$. Any section $\zeta$ of
$\dad{E}$, and in particular any section of $E$, determines a linear function
$\widehat{\zeta}$ on $\ad{E}$ by
$$
\widehat{\zeta}(\varphi)=\pai{\zeta_m}{\varphi}\qquad\text{for every $\varphi\in\ad{E}_m$.}
$$
Then the Poisson bracket is determined by the condition
$$
\{\widehat\zeta_1,\widehat\zeta_2\}=\widehat{[\zeta_1,\zeta_2]},
$$
which for consistency (using linearity and the Leibnitz rule) requires that we
put
$$
\{\widehat\zeta,g\}=\dad{\rho}(\zeta)(g), \qquad\text{and}\qquad \{f,g\}=0,
$$
for $f$ and $g$ functions on $M$.

It is of some interest to mention yet another characterisation of the result
described in Theorem~2. The above Poisson bracket in fact is determined by the
bracket of linear functions coming from sections of $E$, since these span the
set of all linear functions on $\ad{E}$. But the bracket of sections of $E$ is
a section of $V$; it follows that the corresponding Poisson brackets are
independent of the coordinate $\mu_0$, and therefore, $\pd{}{\mu_0}$ is a
symmetry of the Poisson tensor. Conversely, it is obvious that the latter
symmetry property will imply that the bracket of sections in the image of $i$
belongs to the image of $\bs{i}$.

\medskip

It should be noticed that, in general, the Poisson structure we have defined
is not an extension of a Poisson structure on $M\times\R$ by the one on $V^*$.
Indeed, the map $\map{k}{\ad{E}}{V^*}$ is a Poisson map (since it is the dual
of $\bs{i}$ and this is a Lie algebroid morphism), but the map
$\map{l}{M\times\R}{\ad{E}}$ is not a Poisson map, except if $\bs{\rho}=0$.
This is equivalent to the dual mapping $l^*=j$ not being a Lie algebroid
morphism. To prove this, we first show that the brackets
$\{\widehat{\bs{\eta}},f\}\circ l$ and $\{\widehat{\bs{\eta}}\circ l,f\circ
l\}$ are different, except when $\bs{\rho}=0$. Indeed, we have that
$l(m,p)=pe^0$, so that $\widehat{\bs{\eta}}\circ l=0$ and hence
$\{\widehat{\bs{\eta}}\circ l,f\circ l\}=0$. On the other hand,
$\{\widehat{\bs{\eta}},f\}\circ l=\bs{\rho}(\bs{\eta})(f)\circ
l=\bs{\rho}(\bs{\eta})(f)$, since $\bs{\rho}(\bs{\eta})(f)$ is a function on
$M$. Therefore, the two brackets are equal if and only if $\bs{\rho}=0$.
Similarly, one can calculate the brackets $\{\widehat{e_0},f\}\circ
l=\{\mu_0,f\}$ and $\{\widehat{e_0}\circ l,f\circ l\}=\rho(e_0)(f)$. This
imposes that the Poisson tensor on $M\times\R$ has to be
$\Lambda=\pd{}{\mu_0}\wedge X_0$, with $X_0=\rho(e_0)$. (Notice that
$\rho(e_0)=\rho(\sigma)$ for any section $\sigma$ of $E$.) The other brackets
vanish, so there are no further conditions.

\medskip

In coordinates, we have
\begin{align*}
& \{x^i,x^j\}=0 \qquad
&& \{\mu_0,x^i\}=\rho^i_0  \qquad  &&& \{\mu_\alpha,x^i\}=\rho^i_\alpha      \\
& \{\mu_0,\mu_\beta\}=C^\gamma_{0\beta}\mu_\gamma \qquad &&
\{\mu_\alpha,\mu_\beta\}=C^\gamma_{\alpha\beta}\mu_\gamma &&&
\end{align*}
and therefore the Poisson tensor is
\begin{gather*}
\Lambda_{\ad{E}} =\rho^i_\alpha\pd{}{\mu_\alpha}\wedge\pd{}{x^i}+
  \frac{1}{2}\mu_\gamma\,C^\gamma_{\alpha\beta}\pd{}{\mu_\alpha}\wedge \pd{}{\mu_\beta}
  +{}\qquad\qquad\qquad\\
  \qquad\qquad\qquad{}
+\pd{}{\mu_0}\wedge\left(\rho^i_0\pd{}{x^i}
        +\mu_\gamma\,C^\gamma_{0\beta}\pd{}{\mu_\beta}\right).
\end{gather*}

\section{Examples}

\paragraph{The canonical affine Lie algebroid}

The canonical example of a Lie algebroid over an affine bundle is the first
jet bundle $J^1M\rightarrow M$ to a manifold $M$ fibered over the real line
$\map{\prr}{M}{\R}$. The elements of the manifold $J^1M$ are equivalence
classes $j^1_t\gamma$ of sections $\gamma$ of the bundle $\map{\prr}{M}{\R}$,
where two sections are equivalent if they have first
order contact at the point $t$. It is an affine bundle whose associated vector
bundle is $\ver{\prr}$ the set of vectors tangent to $M$ which are vertical
over $\R$. In this case it is well-known that $\ad{J^1M}=T^*M$, and therefore
$\widetilde{J^1M}=TM$. The canonical inmersion is given by
$$
i(j^1_t\gamma)=\dot{\gamma}(t),
$$
i.e.\ it maps the 1-jet of the section $\gamma$ at the point $t$ to the vector
tangent to $\gamma$ at the point $t$. In coordinates, if $j^1_t\gamma$ has
coordinates $(t,x,v)$ then
$$
i(t,x,v)=\pd{}{t}\Big|_{(t,x)}+v^i\pd{}{x^i}\Big|_{(t,x)}
$$
An element $w$ of the associated vector bundle $\ver{\prr}$ is of the form
$$
w=w^i\pd{}{x^i}\Big|_{(t,x)}
$$

The bracket of sections of $J^1M$ is defined precisely by means of the above
identification of 1-jets with vectors $v\in TM$ which projects to the vector
$\partial/\partial t$. In coordinates a section $X$ of $J^1M$ is identified
with the vector field
$$
X=\pd{}{t}+X^i(t,x)\pd{}{x^i},
$$
and the bracket is
$$
[X,Y]=\{X(Y^i)-Y(X^i)\}\pd{}{x^i}
$$
which is obviously a section of the vector bundle.

\paragraph{Affine distributions}
An affine $E$ sub-bundle of $J^1M$ is involutive if the bracket of sections of
the sub-bundle is a section of the associated vector bundle. Therefore, taking
as anchor the natural inclusion into $TM$ and as bracket the restriction of
the bracket in $J^1M$ to  $E$ we have an affine Lie algebroid structure on
$E$.

\paragraph{Lie algebra structures on affine spaces}
We consider the case in which the manifold $M$ reduces to one point $M=\{m\}$.
Thus our affine bundle is $E=\{m\}\times A$  and the associated vector bundle
is $W\equiv\{m\}\times V$ for some affine space $A$ over the vector space $V$.
Then, a Lie algebroid structure over the affine bundle $E$ is just an affine
Lie algebra structure over $A$. Indeed, every section of $E$ and of $W$ is
determined by a point in $A$ and $V$, respectively. The anchor must vanishes
since $TM=\{0_m\}$, so it does not carry any additional information.

\paragraph{Trivial affine algebroids}
By a trivial affine space we mean just a point $A=\{O\}$, and the associated
vector space is the trivial one $V=\{0\}$. The extended affine dual of $A$ is
$A^\dag=\R$ since the only affine maps defined on a space of just a point are
the constant maps. It follows that the extended bi-dual is $\dad{A}=\R$

Given a manifold $M$, we consider the affine bundle $E=M\times\{O\}$ with
associated vector bundle $V=M\times\{0\}$. On $V$ we consider the trivial
bracket $[\,,\,]=0$ and the anchor $\bs{\rho}=0$, and as derivation $D_O$ we
also take $D_O=0$. Now, to construct a Lie algebroid structure on $E$, we take
an arbitrary vector field $X_0$ on $M$ as given and define the map
$\map{\rho}{E}{TM}$ by $\rho(m,O)=X_0(m)$. Then it follows that $\rho$ is
compatible with $D_O$.

The extended dual of $E$ is $E^\dag=M\times\R$ and the extended bi-dual is
$\dad{E}=M\times\R$. We therefore have one section $e^0$ spanning the set of
sections of $E^\dag$, and the dual element $e_0$ (which is just the image
under the canonical immersion of the constant section of value $0$.)

We want to study the associated exterior differential operator and Poisson
bracket.

For the exterior differential operator, since we have 1-dimensional fibre on
$E^\dag$ it follows that $d e^0=0$. On functions $f\in\cinfty{M}$ we have
$df=\rho(e_0)(f)e^0=X_0(f)e^0$.

For the Poisson structure, since the fibre of $E^\dag$ is 1-dimensional, it is
determined by the equation $\{\hat{e_0},f\}=\rho(e_0)(f)$ and
$\hat{e_0}=\mu_0$. We have that the only non-trivial brackets are
$\{\mu_0,f\}=X_0(f)$. Therefore, the Poisson tensor is
$$
\Lambda=\pd{}{\mu_0}\wedge X_0.
$$

\paragraph{Quotient by a group}
If $\map{p}{Q}{M}$ is a principal $G$-bundle and $M$ is fibred over $\R$, then
$E=J^1Q/G\rightarrow M$ is an affine Lie algebroid. The anchor is
$\rho([j^1_t\gamma])=j^1_t(p\circ\gamma)$. The bracket is obtained by
projecting the bracket on $J^1Q$. We have that $\dad{E}$ is the Atiyah
algebroid $=TQ/G$, see for instance~\cite{atiyah.algebroid}.

\paragraph{Affine actions of Lie algebras}

Let $A$ be an affine space endowed with a Lie algebra structure. By an action
of $A$ on a manifold $M$ we mean an affine map
$\map{\phi}{A}{\vectorfields{M}}$, such that
$[\phi(a),\phi(b)]=\bs{\phi}([a,b])$. Then $M\times A\rightarrow M$ has an
affine Lie algebroid structure. The anchor is $\rho(m,\xi)=\phi(\xi)(m)$ and
the bracket can be defined in terms of constant sections: the bracket of two
constant sections $\sigma_i(m)=(m,\xi_i)$, is the constant section
corresponding to the bracket of the values
$$
[\sigma_1,\sigma_2](m)=(m,[\xi_1,\xi_2]_A).
$$

If we consider the Lie algebra $\dad{A}$ then $\dad{A}$ acts also on the
manifold $M$. The extension $\dad{E}$ is the Lie algebroid associated to
the action of $\dad{A}$.

\paragraph{Poisson manifolds with symmetry}

Consider a Poisson manifold $(M,\Lambda)$ and an infinitesimal symmetry
$Y\in\vectorfields{M}$ of $\Lambda$, that is $\Lie{Y}\Lambda=0$ . Take $E$ to
be $T^*M$ with its natural affine structure, where the associated vector
bundle is $V=T^*M$ itself. On $V$ we consider the Lie algebroid structure
defined by the canonical Poisson structure. For a section $\alpha$ of $E$
(i.e.\ a 1-form on $M$) we define the map $\map{D_\alpha}{\sec{V}}{\sec{V}}$
by
$$
D_\alpha\beta=\Lie{Y}\beta+[\alpha,\beta].
$$
Since $Y$ is a symmetry of $\Lambda$, $D_\alpha$ is  a derivation and clearly
satisfies the required compatibility condition. If we further consider the
affine anchor $\map{\rho}{E}{TM}$, determined by
$\rho(\alpha_m)=\Lambda(\alpha_m)+Y_m$, then we have a Lie algebroid structure
on the affine bundle $E$.

In this case, since there is a distinguished section of $E$ (the zero
section), we have that $\ad{E}=TM\times\R$ and $\dad{E}=T^*M\times\R$.

\paragraph{Jets of sections in a groupoid}

Let $G$ be a Lie groupoid over a manifold $M$ with source $\alpha$ and target
$\beta$ (the notation is as in~\cite{atiyah.algebroid}). Let
$T^\alpha_{G^{(0)}} G=\ker T\alpha|_{G^{(0)}}$ be the associated Lie
algebroid, that is, the set of $\alpha$-vertical vectors at points in
$G^{(0)}$ (the set of identities). The anchor is the map $\rho=T\beta$. Assume
that $M$ is further fibred over the real line, $\map{\prr}{M}{\R}$ and
consider the bundle $E=J^\alpha_{G^{(0)}} G$ of 1-jets of sections of
$\prr\circ\beta$ which are $\alpha$-vertical, at points in ${G^{(0)}}$. This
is an affine bundle whose associated vector bundle is $(T^\alpha_{G^{(0)}}
G)^{\mathrm{ver}}$ the set of $(\prr\circ\beta)$-vertical vectors on
$T^\alpha_{G^{(0)}}G$. If $\bs{i}$ is the natural inclusion of
$(T^\alpha_{G^{(0)}} G)^{\mathrm{ver}}$ into $T^\alpha_{G^{(0)}} G$ and we
define the map $\map{j}{T^\alpha_{G^{(0)}} G}{M\times\R}$ by
$j(v)=(\alpha(\tau_G(v)),t(v))$ (where $t=\prr\circ\beta\circ\tau_G$), then we
have the exact sequence of vector bundles over $M$
$$
\seq 0->(T^\alpha_{G^{(0)}} G)^{\mathrm{ver}} -\bs{i}->T^\alpha_{G^{(0)}} G -j-> M\times\R ->0
$$
and $j^{-1}(M\times\{1\})= J^\alpha_{G^{(0)}}G$. Moreover, the bracket of two sections
of $J^\alpha_{G^{(0)}} G$  is vertical over $\R$ from where it follows that
the Lie algebroid structure of $T^\alpha_{G^{(0)}} G$ restricts to a Lie
algebroid structure on the affine bundle $J^\alpha_{G^{(0)}} G$.

\section{Prolongation}

In this section we define the prolongation of a fibre bundle with respect to a
(vector) Lie algebroid. We are primarily interested in the prolongation of the
bundle $E\rightarrow M$, but we will describe explicitly a more general
construction first, since this does not introduce extra complications (see
also \cite{algebraic.category} for generalities).

Let $\map{\mu}{P}{M}$ be a fibre bundle over the manifold $M$ and
$F\rightarrow M$ a Lie algebroid with anchor $\rho$. Using notations
introduced in \cite{LMLA}, we consider the bundle
$\map{\mu_1}{\prol[F]{P}}{P}$ constructed as follows. The manifold
$\prol[F]{P}$ is the total space of $\rho^*(TP)$, the pull-back of $TP$ by
$\rho$, i.e.\
$$
\prol[F]{P}=\set{(z,V)\in F\times TP}{\rho(z)=T\mu(V)}.
$$
The bundle projection we consider, however, is not the usual one
$$
\map{\mu_2}{\prol[F]{P}}{F}\qquad\qquad \mu_2(z,V)=z,
$$
but rather $\map{\mu_1}{\prol[F]{P}}{P}$, defined by
$$
\mu_1(z,V)=\tau_P(V)
$$
where $\tau_P$ is the tangent bundle projection $TP\rightarrow P$.

A section $X$ of $\prol[F]{P}$ is of the form $X(p)=(\theta(p),V(p))$, where
$\theta$ is a section of $F$ along $\mu$ and $V$ is a vector field on $P$. A
section $X$ is projectable if there exists a section $\sigma$ of $F$ such that
$\mu_2\circ X=\sigma\circ\mu$. It follows that $X$ is projectable if and only
if $\theta=\sigma\circ\mu$, and therefore $X$ is of the form
$X(p)=(\sigma(m),A(p))$, with $m=\mu(p)$.

The prolongation of $P$ inherits a Lie algebroid structure from the one on $F$
and the one on $TP$. The anchor is the map
$$
\map{\rho^1}{\prol[F]{P}}{TP}\qquad\qquad \rho^1(z,V)=V,
$$
and the bracket can be defined in terms of projectable sections as follows. If
$X_1$, $X_2$ are two projectable sections given by
$X_k(p)=(\sigma_k(m),V_k(p))$, $k=1,2$, then the bracket $[X_1,X_2]$ is the
section given by
$$
[X_1,X_2](p)=([\sigma_1,\sigma_2](m),[V_1,V_2](p)).
$$
From this expression it is clear that $[X_1,X_2]$ is also a projectable section.

An element $X$ of $\prol[F]P$ is said to be vertical if $\mu_2(X)=0$.
Therefore it is of the form $(0,V)$, for some vertical vector $V$ on $P$. One
should realise, however, that there are sections $X$ of $\prol[F]P$ for which
$\rho^1(X)$ is a vertical vector on $P$, but $X$ itself is not vertical. Such
elements are of the form $(z,V)$, with $z$ in the kernel of $\rho$.

\bigskip

Given a local basis $\{e_a\}$ of sections of $F$ and  local coordinates
$(x^i,u^I)$ on $P$ we can define a basis of sections of $\prol[F]{P}$, as
follows
\begin{align*}
&\baseX_a(p)=\left(e_a(m),\rho^i_a\pd{}{x^i}\Big|_{p}\right)\\
&\baseV_I(p)=\left(0_m,\pd{}{u^I}\right)
\end{align*}
where $m=\mu(p)$. Thus, any element $Z$ of $\prol[F]{P}$ at $p$
$$
Z=\left(z^a e_a(m), (\rho^i_a z^a)\pd{}{x^i}\Big|_{p}
+v^I\pd{}{u^I}\Big|_{p}\right),
$$
can be represented as
$$
Z=z^a\baseX_a+v^I\baseV_{I},
$$
and $(x^i,u^I,z^a,v^I)$ are coordinates  on $\prol[F]{P}$. The brackets of the
elements of this basis are
$$
[\baseX_a,\baseX_b]=C^c_{ab}\baseX_c \qquad [\baseX_a,\baseV_J]=0 \qquad
[\baseV_I,\baseV_J]=0,
$$
and the exterior differential is determined by
$$
dx^i=\rho^i_a\baseX^a, \qquad du^I=\baseV^I,
$$
and
$$
d\baseX^c=-\frac{1}{2}C^c_{ab}\baseX^a\wedge\baseX^b \qquad d\baseV^I=0,
$$
where $\{\baseX^c,\baseV^I\}$ denotes a dual basis.
\bigskip

As a first step towards the situation we are most interested in, we consider
the case where the Lie algebroid $F$ is either $V$ or $\dad{E}$. We will show
that the Lie algebroid structure on $\eprol{P}$ then is precisely of the kind
we are studying in this paper.

\begin{prop}
Let $\map{\bs{I}}{\vprol{P}}{\eprol{P}}$ and $\map{J}{\eprol{P}}{P\times\R}$
be  the maps
$$
\bs{I}(v,W)=(\bs{i}(v),W)\qquad\text{and}\qquad J(z,V)=(\mu_1(z,V),j_m(z)).
$$
Then the following sequence of vector bundles is exact
$$
\seq 0->\vprol{P}-\bs{I}->\eprol{P}-J->P\times\R->0.
$$
\end{prop}
\begin{proof}
We clearly have that $J\circ\bs{I}=0$, so that $\mathrm{Im}(\bs{I})\subset
\mathrm{Ker}(J)$. Moreover, if $(z,V)\in\mathrm{Ker}(J)$ then $j(z)=0$; hence
there exists a $\bs{v}$ such that $z=\bs{i}(\bs{v})$, so that
$(z,V)=(\bs{i}(\bs{v}),V)=\bs{I}(\bs{v},V)$ is in $\mathrm{Im}(\bs{I})$.
\end{proof}

Therefore, the set of points $\jprol{P}=J^{-1}(P\times\{1\})$ is an affine
bundle whose associated vector bundle is $\vprol{P}$. Explicitly, we have
$$
\jprol{P}=\{(i(a),V)\in\eprol{P}\}
   \simeq\set{(a,V)\in E\times TP}{\rho(a)=T\mu(V)}.
$$
Moreover, the Lie algebroid structure of $\eprol{P}$ restricts to $\jprol{P}$,
defining therefore a Lie algebroid structure on that affine bundle. To see
this, we have to prove that the bracket of sections of $\jprol{P}$ is a
section of $\vprol{P}$. For that it is enough to consider projectable
sections, since they form a generating set. If $Z_1$, $Z_2$ are sections of
$\jprol{P}$, projecting to sections $\sigma_1$ and $\sigma_2$ of $E$, then for
every $p\in P$ we have
$$
[Z_1,Z_2](p)=([\sigma_1,\sigma_2](m),[V_1,V_2](p))\qquad\qquad m=\mu(p).
$$
which is an element of $\vprol{P}$, since $[\sigma_1,\sigma_2](m)\in V_m$.

\bigskip

In what follows we further specialise to the case where the bundle $P$ is just
$E$. This is the space where Lagrangian-type systems can be defined, as will
be shown below. The main observation to make in this case is that we have a
canonical map which allows a kind of splitting.

We recall that a splitting of the sequence $\seq
0->V-\bs{i}->\dad{A}-j->\R->0$ is simply a choice of a point of $A$. (Indeed,
if we have a splitting $h$ and we put $h(1)=z$ then $j(z)=j(h(1))=1$, so that
$z=i(a)$ for some $a\in A$). Therefore, once we have fixed a point $a\in A$ we
have two complementary projectors, a `horizontal' one mapping $z$ to
$j(z)i(a)$ and a `vertical' one $\map{\cmap_a}{\dad{A}}{\dad{A}}$ given by
$\cmap_a(z)=z-j(z)i(a)$. Notice that the image of $\cmap_a$ is in $V$, so it
can be considered as a map form $\dad{E}$ to $V$. (Indeed
$j(\cmap_a(z))=j(z-j(z)i(a))=j(z)-j(z)\cdot1=0$, so that $\cmap_a(z)$  is in
the image of $\bs{i}$.)

In the case of an affine bundle, we therefore have a map
$\map{\cmap_a}{\dad{E}_m}{\dad{E}_m}$ for every $a\in E_m$, and thus a map
$\map{\cmap}{\pr^*\dad{E}}{\pr^*\dad{E}}$ given by
$$
\cmap(a,z)=(a,z-j(z)i(a)).
$$
This map is called the canonical map. As mentioned above, we can consider it
as a map $\map{\cmap}{\pr^*\dad{E}}{\pr^*V}$.  In local coordinates, we have
$$
\cmap=(e^\alpha-y^\alpha e^0)\otimes e_\alpha.
$$

An important concept for the study of dynamical systems on affine Lie
algebroids is that of admissible elements. An element $Z\in\eprol{E}$ is said
to be \textsl{admissible} if, roughly, it has the same projection under
$\pr_1$ and $\pr_2$. More precisely, we set
$$
\mathrm{Adm}(E)=\set{Z\in\eprol{E}}{\pr_2(Z)=i(\pr_1(Z))}.
$$
An equivalent characterisation is the following: $Z$ is admissible if and only
if it belongs to $\jprol{E}$ and $\cmap$ vanishes on its projection to
$\pr^*\dad{E}$, which we will denote by $\pr_{12}(Z)$ (cf.\ \cite{LMLA}),
hence
$$
\mathrm{Adm}(E)=\set{Z\in\jprol{E}}{\cmap(\pr_{12}(Z))=0}.
$$
Indeed, if $Z$ is of the form $Z=(z,V_a)$, then the first condition means that
$z$ is in the image of $i$ and the second then further specifies that
$z=i(a)$.

By a contact 1-form we mean a 1-form $\theta$ on $\eprol{E}$ (i.e.\ a
$\cinfty{E}$-linear map from sections of $\pr_1$ to $\cinfty{E}$), which
vanishes on sections whose image lies in $\mathrm{Adm}(E)$. It follows from
the characterisation of admissible elements that contact forms are locally
spanned by
$$
\theta^\alpha=\baseX^\alpha-y^\alpha \baseX^0.
$$
Any 1-form $\theta$ on $\dad{E}$ determines a contact 1-form
$\overline{\theta}$ by means of the canonical map:
$$
\pai{\overline{\theta}}{Z}=\pai{\theta}{\cmap(\pr_{12}(Z))}.
$$
In coordinates, if $\theta$ is of the form $\theta=\theta_0e^0+\theta_\alpha
e^\alpha$ then
$$
\overline{\theta}=\theta_\alpha(\baseX^\alpha-y^\alpha \baseX^0).
$$
Notice that the elements of the basis $\{\theta^\alpha\}$ of contact 1-forms
are of this type: $\theta^\alpha=\overline{e^\alpha}$. We further will need
the affine function $\widehat{\theta}\in\cinfty{E}$ associated to
$i^*(\theta)$. To be precise, there is of course a linear function on
$\dad{E}$ associated to $\theta$, but we will reserve the notation
$\widehat{\theta}$ for its restriction to $E$, meaning that in coordinates:
$$
\widehat{\theta}= \theta_0 + \theta_\alpha y^\alpha.
$$

With these definitions, we can split (the pullback of) a 1-form $\theta$ on
$\dad{E}$ as follows
$$
\pr_2^*\theta=\widehat{\theta}\baseX^0+\overline{\theta}.
$$
This decomposition is important for various calculations in the next section.

\section{Complete and vertical lifts}

The prolongation structure we have specialised to, with $F=\dad{E}$ and $P=E$,
is visualised in the following diagram.

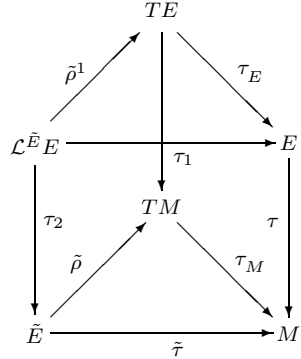
\setlength{\unitlength}{12pt}

\begin{picture}(40,13)(0,2)

\put(7.5,9.5){\vector(1,0){6.6}}  %
\put(11,13){\vector(1,-1){3}}  %
\put(7,10.2){\vector(1,1){2.8}}   %

\put(6.5,9.5){\mybox{$\eprol{E}$}} \put(14.5,9.5){\mybox{$E$}}
\put(10.5,13.7){\mybox{$TE$}}

\put(7,3.5){\vector(1,0){7}}  %
\put(11,7){\vector(1,-1){3}}  %
\put(7,4){\vector(1,1){3}}   %

\put(6.5,3.5){\mybox{$\dad{E}$}} \put(14.5,3.5){\mybox{$M$}}
\put(10.5,7.5){\mybox{$TM$}}

\put(6.5,8.8){\vector(0,-1){4.7}}  %
\put(10.5,13){\vector(0,-1){5}}  %
\put(14.5,9){\vector(0,-1){5}}  %

\put(11,3){\mybox{$\dad{\pr}$}} \put(7.8,5.7){\mybox{$\dad{\rho}$}}
\put(13.3,5.7){\mybox{$\tau^{}_M$}}

\put(11.2,9){\mybox{$\pr_1$}} \put(7.8,11.7){\mybox{$\dad{\rho}^1$}}
\put(13.3,11.7){\mybox{$\tau^{}_E$}}

\put(7.1,7){\mybox{$\pr_2$}} \put(14,7){\mybox{$\pr$}}

\end{picture}

Given a section $\zeta$ of $\dad{E}$ we can define two sections $\zeta^C$ and
$\zeta^V$ of $\eprol{E}$ which are called the complete lift and vertical lift
of $\zeta$, respectively.

The vertical lift is defined by the following sequence of natural
constructions. Given a point $a\in E_m$ and a vector $\bs{v}\in V_m$, we
define the vector $v^V_a\in T_aE$ by its action on functions $f\in\cinfty{E}$:
$$
v^V_a(f)=\frac{d}{dt}f(a+tv)\big\vert_{t=0}.
$$
Next, if $z\in\dad{E}_m$ and $a\in E_m$, the vertical lift of $z$ to the point
$a$ is defined as
$$
z^V_a=(\cmap_a(z))^V_a.
$$
For the special case that $z=\bs{i}(\bs{v})$, this is consistent with the
preceding step:  $z^V_a=\bs{v}^V_a$. Let now
$\map{\xi^V}{\pr^*\dad{E}}{\eprol{E}}$ be the \textsl{vertical lift} map,
determined by $\xi^V(a,z)=(0_{\pr(a)},z^V_a)$. The final step in the
construction now is obvious: if $\zeta$ is a section of $\dad{E}$, we define
the section $\zeta^V$ of $\vprol{E}$, called the vertical lift of $\zeta$, by
putting
$$ \zeta^V(a)=\xi^V(a,\zeta(m)), \qquad\qquad\text{with $a\in E$ and
$m=\pr(a)$.}
$$
In coordinates, if $\zeta=\zeta^0 e_0+\zeta^\alpha e_\alpha$ then
$$
\zeta^V=(\zeta^\alpha-y^\alpha\zeta^0)\baseV_\alpha.
$$

Much of the structure here discussed is encoded in the definition of a
vertical endomorphism $S$ on $\eprol{E}$ (or sections of it), which goes as
follows: $S=\xi^V\circ\cmap\circ\pr_{12}$. Explicitly, if $Z=(z,V_a)$, then $
S(z,V_a)=\xi^V(\cmap_a(z))$. In coordinates, the type (1,1) tensor field $S$
reads
$$
S=(\baseX^\alpha-y^\alpha\baseX^0)\otimes\baseV_\alpha.
$$

The complete lift of a section $\zeta$ of $\dad{E}$ is defined by the following
two conditions which completely characterise it:
\begin{itemize}
\item $\zeta^C$ projects to $\zeta$, i.e.\
$\pr_2\circ\zeta^C=\zeta\circ\pr$, and
\item $\zeta^C$ preserves the set of contact forms, that is, if $\theta$ on $\eprol{E}$ is a
contact form then $d_{\zeta^C}\theta=[i_{\zeta^C},d]\theta$ is contact.
\end{itemize}

In the case of the pullback of a 1-form on $\dad{E}$, making use of the
decomposition as sum of a contact plus a non-contact form
$\pr_2^*\theta=\widehat{\theta}\baseX^0+\overline{\theta}$, and taking into
account that $\zeta^C$ projects onto $\zeta$ and that $\pr_2$ is a morphism of
Lie algebroids, one can verify that
\begin{align*}
d_{\zeta^C}\overline{\theta}
&=\overline{d_{\zeta}\theta}+\widehat{\theta}\,\,\overline{d_\zeta e^0},\\
d_{\zeta^C}\widehat{\theta}
&=\widehat{d_{\zeta}\theta}-\widehat{\theta}\,\widehat{d_\zeta e^0}.
\end{align*}
In fact, any of these two conditions is equivalent to the second condition in
our definition of complete lift.

The coordinate expression of the complete lift of the section
$\zeta=\zeta^0e_0+\zeta^\alpha e_\alpha$ is
$$
\zeta^C=\zeta^0\baseX_0+\zeta^\alpha \baseX_\alpha +
[(\dot{\zeta}^\alpha-y^\alpha\dot{\zeta}^0) +
C^\alpha_\beta(\zeta^\beta-y^\beta\zeta^0)] \baseV_\alpha,
$$
where $C^\alpha_\beta=C^\alpha_{0\beta}+C^\alpha_{\gamma\beta}y^\gamma$ and,
for a function $f\in\cinfty{M}$, the complete lift $\dot{f}\in\cinfty{E}$ is
defined by $\dot{f}=\widehat{df}$. The first two terms of $\zeta^C$ are
determined by the projectability condition, whereas the third term can be
obtained by applying the preceding formula to $\theta=e_\alpha$.

\medskip

The vertical and complete lift satisfy the following properties
\begin{align*}
&d_{\zeta^V}f=0
&&d_{\zeta^V}\widehat\theta=i_{\zeta^C}\bar{\theta}\\
&d_{\zeta^C}f=d_\zeta f
&&d_{\zeta^C}\widehat{\theta}
=\widehat{d_{\zeta}\theta}-\widehat{\theta}\,\widehat{d_\zeta e^0}
\end{align*}
for $f\in\cinfty{M}$ and $\theta$ a 1-form on $\dad{E}$. We prove only
the third; if $\bs{v}=\cmap_a(\zeta)$ then
\begin{align*}
d_{\zeta^V}\widehat\theta(a)
&=\zeta(m)^V_a\widehat\theta
 =\frac{d}{dt}\widehat\theta(a+t\bs{v})\big\vert_{t=0}
 =\frac{d}{dt}(\widehat\theta(a)+ t\pai{\bs{\theta}_m}{\bs{v}})\big\vert_{t=0}\\
&=\pai{\bs{\theta}_m}{\bs{v}}
 =\pai{\bs{\theta}_m}{\cmap_a(\zeta)}
 =(i_{\zeta^C}\bar{\theta})(a).
\end{align*}

Also, it follows form the definition of the vertical endomorphism that
$$
S(\zeta^C)=\zeta^V\qquad\text{\and}\qquad S(\zeta^V)=0.
$$

Using the above equations it is a matter of a routine calculation to prove the
following commutation relations.
\begin{align*}
&[\zeta_1^C,\zeta_2^C]=[\zeta_1,\zeta_2]^C\\
&[\zeta_1^C,\zeta_2^V]=[\zeta_1,\zeta_2]^V+
  {\pai{\zeta_1}{e^0}}\dot{}\,\zeta_2^V\\
&[\zeta_1^V,\zeta_2^V]=\pai{\zeta_1}{e^0}\zeta_2^V-\pai{\zeta_2}{e^0}\zeta_1^V.
\end{align*}

The above definitions and relations are greatly simplified if we restrict to
sections of the associated vector bundle $V$. Indeed, if
$\bs{\sigma}=\sigma^\alpha e_\alpha$ is a section of $V$, then the complete
and vertical lifts of $\bs{\sigma}$ have the coordinate expressions
$$
\bs{\sigma}^C=\sigma^\alpha \baseX_\alpha +
(\dot{\sigma}^\alpha + C^\alpha_\beta\sigma^\beta) \baseV_\alpha
\qquad\qquad
\bs{\sigma}^V=\sigma^\alpha \baseV_\alpha,
$$
and the action of the complete lift over linear functions and contact forms is
given by
$$
d_{\bs{\sigma}^C}\overline{\theta}
=\overline{d_{\bs{\sigma}}\theta}
\qquad\qquad
d_{\zeta^C}\widehat{\theta}
=\widehat{d_{\bs{\sigma}}\theta},
$$
since $d_{\bs{\sigma}} e^0=0$. Furthermore, the commutation relations are as in
the usual vector Lie algebroid case:
$$
[\bs{\sigma}^C,\bs{\eta}^C]=[\bs{\sigma},\bs{\eta}]^C\qquad
[\bs{\sigma}^C,\bs{\eta}^V]=[\bs{\sigma},\bs{\eta}]^V\qquad
[\bs{\sigma}^V,\bs{\eta}^V]=0.
$$

\section{Lagrangian-type systems on an affine algebroid}

In this section, we first consider dynamical systems on $E$, which
geometrically are defined in the way standard second-order differential
equations on a tangent bundle or first-jet bundle are conceived, but do not
necessarily correspond, locally, to second-order equations. As in
\cite{affine.algebroid} therefore, we will call them \sode s. We shall
subsequently discuss a class of \sode s, which come from a (constrained)
variational problem and therefore are said to be of Lagrangian type.

A curve $\map{\gamma}{\R}{E}$ is said to be admissible if
$$
\rho\circ\gamma=\dot{\gamma}_{\scriptscriptstyle M},
$$
where $\gamma_M=\pr\circ\gamma$ is the projected curve to the base. In
coordinates, if $\gamma(t)=(x^i(t),y^\alpha(t))$ then $\gamma$ is admissible if
$$
\dot{x}^i(t)=\rho^i_0(x(t))+\rho^i_\alpha(x(t))y^\alpha(t).
$$
A curve is admissible if and only if  its prolongation takes values in the set
of admissible elements. Indeed, the prolongation of the curve $\gamma$ is the
curve $\gamma^c(t)=(i(\gamma(t)),\dot{\gamma}(t))$, and this is an admissible
element if and only if it is in $\eprol{E}$, that is
$\rho(\gamma(t))=T\pr(\dot{\gamma})(t)=\dot{\gamma}_{\scriptscriptstyle
M}(t)$.

\begin{dfn}
A \sode\ on $E$ is a section $\Gamma$ of $\mathrm{Adm}(E)$, that is, a section
of $\eprol{E}$ which takes values in the set of admissible elements.
\end{dfn}

From this definition, it readily follows that $\pai{\Gamma}{\baseX^0}=1$ and
the integral curves of $\dad{\rho}^1(\Gamma)$ are admissible curves.
Conversely, any section $Z$ of $\eprol{E}$ such that the integral curves of
$\dad{\rho}^1(Z)$ are admissible is a  \sode. From  the alternative
characterisation of $\mathrm{Adm}(E)$ as the set of elements of $\jprol{E}$
which vanish under $\cmap$, it follows that a section $\Gamma$ of $\eprol{E}$ is
a \sode\ if and only if $S(\Gamma)=0$ and $\pai{\Gamma}{\baseX^0}=1$.

Locally, a \sode\ $\Gamma$ is of the form
$$
\Gamma=\baseX_0+y^\alpha\baseX_\alpha+F^\alpha\baseV_\alpha.
$$
and the vector field $\dad{\rho}^1(\Gamma)$ is of the form
$$
\dad{\rho}^1(\Gamma) =(\rho^i_0+\rho^i_\alpha
y^\alpha)\pd{}{x^i}+F^\alpha\pd{}{y^\alpha}.
$$

\bigskip
Now, to define Lagrangian-type equations in a coordinate free way, we can (as
in \cite{LMLA}) simply mimic the usual construction on a first-jet bundle. For
a given function $L$ on $E$, we define the Cartan 1-form $\Theta_L$ on
$\eprol{E}$ by
$$
\Theta_L=dL\circ S+L\baseX^0
$$
and the Cartan 2-form $\Omega_L$ by $\Omega_L=-d\Theta_L$. We say that a
\sode\ $\Gamma$ is of Lagrangian type if
$$
i_\Gamma\Omega_L=0.
$$
If the Lagrangian is regular (the 2-form $\Omega_L$ has maximal rank at every
point) then there are no other solutions than multiples of a \sode, but in the
singular case this is a condition to be imposed.

In coordinates, we get
$$
\Theta_L=\pd{L}{y^\alpha}\theta^\alpha+L\baseX^0,
$$
and the expression of $\Omega_L$ is simplified by fixing an arbitrary \sode\
$\Gamma_0=\baseX_0+y^\alpha\baseX_\alpha+F_0^\alpha\baseV_\alpha$ and using
the basis
$\{\baseX^0,\theta^\alpha,\psi^\alpha=\baseV^\alpha-F_0^\alpha\baseX^0\}$ dual
to $\{\Gamma_0,\baseX_\alpha,\baseV_\alpha\}$,
\begin{align*}
\Omega_L &=\left(
     d_{\Gamma_0}\left(\pd{L}{y^\alpha}\right)-\pd{L}{y^\gamma}C^\gamma_\alpha
        -\rho^i_\alpha\pd{L}{x^i} \right)\theta^\alpha\wedge\baseX^0 \\
&{}+\pd{^2L}{y^\alpha\partial y^\beta}\theta^\alpha\wedge\psi^\beta
+\frac{1}{2}\left(
     \rho^i_\beta\pd{^2L}{x^i\partial y^\alpha}
    -\rho^i_\alpha\pd{^2L}{x^i\partial y^\beta}
    +\pd{L}{y^\gamma}C^\gamma_{\alpha\beta}
\right)\theta^\alpha\wedge\theta^\beta.
\end{align*}
Lagrangian-type \sode-equations are of the form
\begin{align*}
\dot x^i &= \rho^i_0 + \rho^i_\alpha y^\alpha, \\
\frac{d}{dt}\Big(\pd{L}{y^\alpha}\Big) &= \rho^i_\alpha \pd{L}{x^i} +
C^\gamma_{\alpha}\pd{L}{y^\gamma}.
\end{align*}

It is interesting to verify that such equations also can be obtained from a
geometric calculus of variations approach. We explain how this works without
working out all the technical details. Given a function $L\in\cinfty{E}$ and
two points $m_0$ and $m_1$ on $M$, consider the problem of determining the
critical curves of the functional
$$
J(\gamma)=\int_\gamma L
=\int_{t_0}^{t_1}L(\gamma(t))\,dt
$$
defined on the set of admissible curves $\map{\gamma}{[t_0,t_1]}{E}$, for
which $\gamma_{\scriptscriptstyle M}$ in the base manifold has fixed endpoints
$m_0$ and $m_1$. This is a constrained problem, since the curves we consider
are restricted to be admissible, i.e.\ they have to satisfy the constraints
$\dot{x}^i=\rho^i_0+\rho^i_\alpha\, y^\alpha$. We should therefore be more
specific about the class of admissible variations we will allow; they will be
generated by complete lifts of sections of $V$, as follows.

Let $\bs{\sigma}$ be a section of $V$ such that
$\bs{\sigma}(m_0)=\bs{\sigma}(m_1)=0$. We consider the vector fields
$X=\bs{\rho}(\bs{\sigma})$ and $Y=\bs{\rho}^1(\bs{\sigma}^C)$, and we denote
their flows by $\psi_s$ and $\Psi_s$, respectively. It follows that
$\psi_s(m_0)=m_0$ and $\psi_s(m_1)=m_1$. The family of curves $\chi(s,t)=
\Psi_s(\gamma(t))$ is a 1-parameter family of admissible variations of
$\gamma$: that $\chi(s,t)$ projects onto $\chi_{\scriptscriptstyle M}(s,t)=
\psi_s(\gamma_{\scriptscriptstyle M}(t))$ is obvious; the fact that
$\chi(s,t)$ is an admissible curve for every fixed $s$ requires more work and
is left to the reader. At the endpoints $t_0$ and $t_1$, we have
$$
\chi(s,t_i)=\Psi_s(\gamma(t_i))=\psi_s(m_i)=m_i.
$$

The infinitesimal variation fields we consider are of the form
$Z=Y\circ\gamma$; their projection to $M$ is
$W=X\circ\gamma_{\scriptscriptstyle M}$. Therefore, the variation of $L$ along
$\chi(s,t)$ at $s=0$ is given by
$$
\pd{(L\circ\chi)}{s}(0,t)=Z(t)(L)=Y(L)(\gamma(t))
=\rho^1(\bs{\sigma}^C)(L)(\gamma(t))=d_{\bs{\sigma}^C}L(\gamma(t)),
$$
from which it follows that
$$
\frac{d}{ds}J(\chi_s)\Big|_{s=0}
=\int_{\gamma}d_{\bs{\sigma}^C}L.
$$
If $\bs{\sigma}$ is a section satisfying the conditions given above, then so is
$f\bs{\sigma}$ for every function $f$ on $M$. Taking into account that
$(f\bs{\sigma})^C=f\bs{\sigma}^C+\dot{f}\bs{\sigma}^V$, we have that
\begin{align*}
0
&=\int_{\gamma}d_{(f\bs{\sigma})^C}L\\
&=\int_{\gamma}fd_{\bs{\sigma}^C}L+\dot{f}d_{\bs{\sigma}^V}L\\
&=f\pai{dL}{\bs{\sigma}^V}\Big|_{\gamma(t_0)}^{\gamma(t_1)} +\int_{\gamma}f\{
d_{\bs{\sigma}^C}L -d_{\Gamma}\pai{dL\circ S}{\bs{\sigma}^C}\}\\
&=\int_{\gamma}fi_{\bs{\sigma}^C}\{dL-d_{\Gamma}(dL\circ S)\},
\end{align*}
whereby $\Gamma$ is the \sode\ of which the extremals we are looking for will
be solutions, and we have made use of the property that $[\sigma^C,\Gamma]$ is
vertical, as one can easily verify in coordinates.

Since $f$ is arbitrary, the fundamental lemma of the calculus of variations
implies that its coefficient must vanish along extremals $\gamma(t)$ and
therefore also in an open neighbourhood on $E$. So the vanishing of the
variation of $J$ is equivalent to
$$
i_{\bs{\sigma}^C}(dL-d_{\Gamma}(dL\circ S))=0.
$$
One easily verifies that $d_{\Gamma}(dL\circ S)-dL$ is `semi-basic', and since
it vanishes on the complete lift of arbitrary sections $\bs{\sigma}$ of $V$,
it follows that $d_{\Gamma}(dL\circ S)-dL=\lambda \baseX^0$. The value of
$\lambda$ can be found by contraction with $\Gamma$:
$$
\lambda
=\pai{d_{\Gamma}(dL\circ S)-dL}{\Gamma}
=i_\Gamma d_{\Gamma}(dL\circ S)-d_\Gamma L
=d_{\Gamma}i_\Gamma(dL\circ S)-d_\Gamma L
=-d_\Gamma L.
$$
Thus the Euler-Lagrange equations can be written as $d_\Gamma\Theta_L=dL$,
from which it follows, since $di_\Gamma\Theta_L=dL$, that
$i_\Gamma\Omega_L=0$.

To finish, we also outline briefly the relation to a Hamilton-type
formulation. For that we consider the prolongation of the extended dual, i.e.\
$\eprol{\ad{E}}$. In this bundle we have a canonical 1-form $\theta$
(cf.~\cite{medina} or~\cite{HMLA}) given by
$$
\pai{\theta_0}{(z,V_\varphi)}=\pai{z}{\varphi}.
$$
In a local basis $\{\baseX^0,\baseX^\alpha,\baseP^0,\baseP^\alpha\}$ of
$\eprol{\ad{E}}$ induced by a frame on $E$, the 1-form $\theta_0$ reads
$$
\theta_0=\mu_0\baseX^0+\mu_\alpha\baseX^\alpha.
$$
The canonical symplectic form $\omega_0$ is the 2-form $\omega_0=-d\theta_0$,
which in coordinates has the expression
$$
\omega_0=\baseX^0\wedge\baseP_0+\baseX^\alpha\wedge\baseP_\alpha
+\mu_\gamma C^\gamma_{0\beta}\baseX^0\wedge\baseX^\beta
+\frac{1}{2}\mu_\gamma C^\gamma_{\alpha\beta}\baseX^\alpha\wedge\baseX^\beta.
$$

The Lagrangian defines a map $\mathcal{F}_L$ from $E$ to $\ad{E}$, the
Legendre transformation, defined as follows. The element
$\mathcal{F}_L(a)\in\ad{E}$ is the affine approximation  (first order Taylor
polynomial) of $L$ at $a$. In other words, if $v$ is the vector such that
$b=a+\bs{v}$ and $g$ is the function $g(t)=L(a+t\bs{v})$ then
$$
\mathcal{F}_L(a)(b)=g(0)+g'(0).
$$
In coordinates
$$
\mathcal{F}_L(x^i,y^\alpha)
=(x^i,L-\pd{L}{y^\alpha}y^\alpha,\pd{L}{y^\alpha})
$$
Then we have, as in the standard theory, that the Cartan forms are the the
pullback of the canonical forms by the prolongation of the Legendre
transformation:
$$
(\prol{\mathcal{F}_L})^*\theta_0=\Theta_L \qquad\text{and}\qquad
(\prol{\mathcal{F}_L})^*\omega_0=\Omega_L.
$$
where $\map{\prol{\mathcal{F}_L}}{\eprol{E}}{\eprol{\ad{E}}}$ is the map
$$
\prol{\mathcal{F}_L}(z,V)=(z,T\mathcal{F}_L(V)).
$$

\section{Final comments}

The case we studied in~\cite{affine.algebroid} is when $M$ is fibred over the
real line $\map{\prr}{M}{\R}$ and the image of the anchor map $\rho$ belongs
to $i(J^1M)\subset TM$. As said before, the motivation to investigate this
particular case comes from potentially interesting applications of a
time-dependent generalisation of the by now classical `Lagrangian systems' on
(vector) Lie algebroids (cf.\ \cite{Wein}). Needless to say, we should be able
to recover this case simply within the present more general set-up.
Essentially, in the case of the extra fibration, we have coordinates
$(t,x^i)$, and we can think of $t$ as being the zeroth coordinate. Then, we
have $\rho^0_0=1$ and $\rho^0_i=0$, so that
$$
\rho(e_0+y^\alpha e_\alpha)=\pd{}{t}+(\rho^i_0+\rho^i_\alpha
y^\alpha)\pd{}{x^i}.
$$

An interesting feature is that in this case we have $e^0=dt$, so that $e^0$ is
not only closed but also exact. In fact, as argued in Section~6, this is the
condition for a Lie algebroid structure on an affine bundle to have a
1-jet-valued anchor. There are corresponding changes, for example in the
formula for the exterior derivative of a function and in the fundamental
Poisson brackets on $\ad{E}$. We refer to \cite{affine.algebroid} for this.

For the sake of clarity, however, it is useful to point out a rather subtle
difference between the first and the present approach. In
\cite{affine.algebroid}, we also discussed the notion of prolongation.
Compared to the diagram we have here, in Section~9, the difference is that the
bottom line was $E\rightarrow M$ rather than $\dad{E}\rightarrow M$ and the
prolonged bundle accordingly was an affine bundle rather than the vector
bundle $\eprol{E}\rightarrow E$ we have here. In fact, for the purpose of
defining (geometrically) time-dependent systems of Lagrangian type, one can
construct all the necessary tools also on the affine prolongation of
\cite{affine.algebroid}. However, we wish to discuss in a forthcoming paper
aspects such as the non-linear and linear connections which are naturally
associated to a \sode\ on an affine bundle $E$, and for that purpose, even
specifically for the time-dependent framework, it turns out to be much more
appropriate to use the prolongation structure of the diagram of Section~9.

\end{document}